\documentclass[12pt,a4paper,final]{amsart}
 
\usepackage[latin1]{inputenc} 
 
\usepackage[notcite,notref,color]{showkeys} 

\usepackage{amsmath,amssymb,amsthm,bbm} 
\usepackage{enumitem} 
\usepackage[overload]{textcase}
 
\usepackage{epsfig,psfrag,color} 

\usepackage[colorlinks,bookmarks,linkcolor=black,citecolor=black]{hyperref} 
\usepackage{verbatim} 
 
\setlist{itemsep=1pt,parsep=0pt,topsep=2pt,partopsep=0pt}  
\setenumerate{leftmargin=*} 
 
\def\itm#1{\rm ({#1})} 
\def\itmit#1{\itm{\it #1\,}} 
\def\rom{\itmit{\roman{*}}} 
\def\abc{\itmit{\alph{*}}}

\newcommand{\cF}{\mathcal{F}}
\newcommand{\cFs}{\mathcal{F}^{\stern}}

\let\subset\subseteq  
 
\let\rho\varrho 
\let\theta\vartheta

\def\dcup{\mathbin{\dot{\cup}}}  

\def\sg{\sqrt{\gamma}}

\let\polishlcross=\l
\def\l{\ifmmode\ell\else\polishlcross\fi}

\def\qqand{\qquad\text{and}\qquad}

\newtheorem{theorem}{Theorem}
\newtheorem{lemma}[theorem] {Lemma}    
    
\newtheorem{conjecture}[theorem] {Conjecture}

\newtheorem{fact}[theorem]{Fact}

\theoremstyle{remark}

\newcommand{\By}[2]{\overset{\mbox{\tiny{#1}}}{#2}} 
\newcommand{\ByRef}[2]{   \By{\eqref{#1}}{#2} }

\newcommand{\leByRef}[1]{ \ByRef{#1}{\le} } 
\newcommand{\geByRef}[1]{ \ByRef{#1}{\ge} }

\newcommand{\EMAIL}[1]{  \textit{E-mail}: \texttt{#1} } 

\newcommand{\bigdcup}{\mathbin{\text{\mbox{\makebox[0mm][c]{\hphantom{$\bigcup$}$\cdot$}$\bigcup$}}}}
\newcommand{\stern}{\vphantom{a}^*}

\DeclareMathOperator{\ex}{ex}
\DeclareMathOperator{\biex}{biex}

\title[An improved error term for minimum $H$-decompositions]{
  An improved error term for minimum $H$-decompositions of graphs
}

  \author[Peter Allen]{Peter Allen*}
  \author[Julia B\"ottcher]{Julia B\"ottcher*}

  \thanks{
    *
    Instituto de Matem\'atica e Estat\'{\i}stica, Universidade de
    S\~ao Paulo, Rua do Mat\~ao 1010, 05508--090~S\~ao Paulo, Brazil.
   \EMAIL{allen|julia@ime.usp.br}
 }

  \author[Yury Person]{Yury Person\dag}
  \thanks{
    \dag\ 
     Institut f\"ur Mathematik, Freie Universit\"at Berlin, Arnimallee 3, 14195 Berlin, Germany
   \EMAIL{person@math.fu-berlin.de}
  }

  \thanks{
    PA was partially supported by FAPESP (Proc.~2010/09555-7);
    JB by FAPESP (Proc.~2009/17831-7);
    PA and JB by CNPq (Proc.~484154/2010-9);
    YP by GIF grant no.~I-889-182.6/2005.
    The cooperation of the three authors was
    supported by a joint CAPES-DAAD project (415/ppp-probral/po/D08/11629,
    Proj.~no.~333/09).
    The authors are grateful to NUMEC/USP, N\'ucleo de Modelagem Estoc\'astica e
    Complexidade of the University of S\~ao Paulo, and Project MaCLinC/USP, for supporting this
    research.
  } 
 
\date{\today} 

\begin{document} 

\begin{abstract} 
  We consider partitions of the edge set of a graph~$G$ into copies of a fixed
  graph~$H$ and single edges.  Let $\phi_H(n)$ denote the minimum
  number~$p$ such that any $n$-vertex~$G$ admits such a partition with at
  most~$p$ parts.  We show that
  $\phi_H(n)=\ex(n,K_r)+\Theta\big(\biex(n,H)\big)$ for $\chi(H)\ge 3$, where $\biex(n,H)$ is the extremal number of the
  decomposition family of $H$. Since $\biex(n,H)=O(n^{2-\gamma})$ for some
  $\gamma>0$ this improves on the bound $\phi_H(n)=\ex(n,H)+o(n^2)$ by
  Pikhurko and Sousa [J. Combin. Theory Ser. B 97 (2007), 1041--1055]. In
  addition it extends a result of \"Ozkahya and Person [J. Combin. Theory
  Ser. B, to appear].
\end{abstract} 

\maketitle

\setcounter{footnote}{0}
\renewcommand{\thefootnote}{\fnsymbol{footnote}}
\section{Introduction} \label{sec:Intro}

We study edge decompositions of a graph~$G$ into disjoint copies of another
graph~$H$ and single edges.  More formally, an
\emph{$H$-de\-com\-po\-si\-tion} of~$G$ is a decomposition
$E(G)=\bigdcup_{i\in[t]} E(G_i)$ of its edge set, such that for
all~$i\in[t]$ either $|E(G_i)|=1$ or~$G_i$ is isomorphic to~$H$. Let
$\phi_H(G)$ denote the minimum~$t$ such there is a decomposition
$E(G)=\bigdcup_{i\in[t]} E(G_i)$ of this form, and let $\phi_H(n) :=
\max_{v(G)=n} \phi_H(G)$.


The function $\phi_H(n)$ was first studied in the seventies by Erd\H{o}s, Goodman and
P\'osa~\cite{EGP66}, who showed that the minimal number~$k(n)$ such that
every $n$-vertex graph admits an edge decomposition into~$k(n)$ cliques
equals $\phi_{K_3}(n)$. They also proved that $\phi_{K_3}(n)=\ex(n,K_3)$,
where $\ex(n,H)$ is the maximum number of edges in an $H$-free graph on $n$
vertices. A decade later this result was extended to~$K_r$ for arbitrary~$r$ by
Bollob\'as~\cite{Bo76} who showed that $\phi_{K_r}(n)=\ex(n,K_r)$ for all
$n\ge r\ge 3$.

General graphs~$H$ were considered only recently by Pikhurko and
Sousa~\cite{PS07}, who proved the following upper bound for $\phi_H(n)$.

\begin{theorem}[Theorem~$1.1$ from~\cite{PS07}]\label{thm:PSthm}
  If $\chi(H)=r\ge 3$ then
 \[
  \phi_H(n)=\ex(n,K_r)+o(n^2).
  \]
\end{theorem}

Pikhurko and Sousa also conjectured that if
$\chi(H)\ge3$ and if~$n$ is sufficiently large, then the correct value is the
extremal number of $H$.

\begin{conjecture}\label{conj1}
  For any graph $H$ with chromatic number at least $3$, there is an
  $n_0=n_0(H)$ such that $\phi_H(n)=\ex(n,H)$ for all $n\ge n_0$.
\end{conjecture}

We remark that the function $\ex(n,H)$ is known precisely only for some
graphs~$H$, which renders Conjecture~\ref{conj1} difficult.  However,
$\ex(n,H)$ is
known for the family of \emph{edge-critical} graphs~$H$, that is, graphs
with $\chi(H)>\chi(H-e)$ for some edge~$e$.  And in fact, after
Sousa~\cite{sousa2010+,sousa2011,sousa2005} proved Conjecture~\ref{conj1}
for a few special edge-critical graphs, \"Ozkahya and Person~\cite{OP11}
verified it for all of them.

Our contribution is an extension of the result of \"Ozkahya and Person
to arbitrary graphs~$H$, which also improves on Theorem~\ref{thm:PSthm}. We
need the following definition.
Given a graph~$H$ with $\chi(H)=r$, the \emph{decomposition family} $\cF_H$
of $H$ is the set of bipartite graphs which are obtained from $H$ by
deleting $r-2$ colour classes in some $r$-colouring of~$H$. Observe
that~$\cF_H$ may contain graphs which are disconnected, or even have
isolated vertices.  Let~$\cF^{\stern}_H$ be a minimal subfamily of~$\cF_H$
such that for any $F\in\cF_H$, there exists $F'\in\cF^{\stern}_H$ with
$F'\subseteq F$. We define
\begin{equation*}
  \biex(n,H):=\ex(n,\cF_H)=\ex(n,\cF^{\stern}_H) \,.
\end{equation*}
Our main result states that the $o(n^2)$ error term in
Theorem~\ref{thm:PSthm} can be replaced by $O\big(\biex(n,H)\big)$, which is
$O(n^{2-\gamma})$ for some $\gamma>0$ by the result of K\"ovari, Tur\'an
and S\'os~\cite{KovSosTur54}. Furthermore, we show that our error term is of the
correct order of magnitude.

\begin{theorem}
\label{thm:biex}
  For every integer~$r\ge3$ and every graph~$H$ with $\chi(H)=r$ there are
  constants~$c=c(H)>0$ and $C=C(H)$ and an integer~$n_0$ such that for all $n\ge
  n_0$ we have
  \begin{equation*}
    \ex(n,K_r)+c\cdot\biex(n,H)\le \phi_H(n) \le \ex(n,K_r)+C\cdot\biex(n,H) \,.
  \end{equation*}
\end{theorem}

Since for every edge-critical~$H$ and every~$n$ we have $\biex(n,H)=0$,
this is indeed an extension of the result of \"Ozkahya and Person. 

\section{Outline of the proof and auxiliary lemmas}

The lower bound of Theorem~\ref{thm:biex} is obtained as follows. We let $F$ be
an $n$-vertex $\cF^{\stern}_H$-free graph with $\biex(n,H)$ edges, and let
$c=(r-1)^{-2}$. There is an $n/(r-1)$-vertex subgraph $F'$ of $F$ with at least
$c\cdot e(F)$ edges. We let $G$ be obtained from the complete balanced
$(r-1)$-partite graph on $n$ vertices by inserting $F'$ into the largest part.
Clearly, we have $e(G)\ge\ex(n,K_r)+c\cdot\biex(n,H)$, and by definition of $\cF^{\stern}_H$, the
graph $G$ is $H$-free, and therefore satisfies
$\phi_H(G)=e(G)\ge\ex(n,K_r)+c\cdot\biex(n,H)$.

The upper bound of Theorem~\ref{thm:biex} is an immediate consequence of the
following result.

\begin{theorem}\label{thm:main}
  For every integer~$r\ge 3$ and every graph~$H$ with $\chi(H)=r$ there is a
  constant~$C=C(H)$ and an integer~$n_0$ such that the following holds.
  Every graph~$G$ on $n\ge n_0$ vertices and with
  \begin{equation*}
    e(G)\ge\ex(n,K_r)+C\cdot\biex(n,H)
  \end{equation*}
  satisfies $\phi_H(G)\le\ex(n,K_r)$.
\end{theorem}

The proof of this theorem (see Section~\ref{sec:proof}) uses the
auxiliary lemmas collected in this section and roughly proceeds as follows.
We start with a graph $G=(V,E)$ on $n$ vertices with $e(G)\ge
\ex(n,K_r)+C\cdot\biex(n,H)$.
For contradiction we assume that
$\phi_H(G)>\ex(n,K_r)$.
This allows us to use a stability-type result (Lemma~\ref{lem:stab}), which
supplies us with a partition $V=V_1,\ldots, V_{r-1}$ with parts of roughly
the same size and with few edges inside each part. Since $e(G)\ge
\ex(n,K_r)+C\cdot\biex(n,H)$ we also know that between two parts only few
edges are missing.
Next, in each part~$V_i$ we identify the (small) set~$X_i$ of those vertices with
many edges to~$V_i$ and set $V'_i:=V_i\setminus X_i$ and $X:=\bigcup
X_i$. 

Then we consider the graph $G[V\setminus X]$, and identify a copy of
some~$F$ in the decomposition family of~$H$ in any $G[V'_i]$, which we then
complete to a copy of~$H$ using the classes~$V'_j$ with $j\neq i$ (see
Lemma~\ref{lem:lowdegree}). We delete this copy of~$H$ from~$G$ and repeat
this process. We shall show that this is possible until the number of edges in all
$V_i\setminus X_i$ drops below $\biex(n,H)$, and thus this gives many
edge-disjoint copies of~$H$ in~$G$.

Finally, we find edge-disjoint copies of~$H$ each of which has one of
its colour classes in~$X$ and the other $(r-1)$ colour classes in $V'_1$,
\dots, $V'_{r-1}$ (see Lemma~\ref{lem:highdegree}). It is possible to find
many copies of~$H$ in this way, because every vertex in~$X$ has many
neighbours in \emph{every}~$V_i$.

In total these steps will allow us to find enough $H$-copies to obtain a
contradiction.

\subsection{Notation}

Let~$G$ be a graph and $V(G)=V_1\dcup\dots\dcup V_s$ a partition of its
vertex set. We write $e(V_i)$ for the number of edges of~$G$ with both ends
in~$V_i$ and $e(V_i,V_j)$ for the number of edges of~$G$ with one end
in~$V_i$ and one end in~$V_j$. Moreover, for~$v\in V$ we let
$\deg_G(v,V_i)=\deg(v,V_i)$ denote the number of neighbours of~$v$
in~$V_i$.
An edge of~$G$ is called \emph{crossing} (for
$V_1\dcup\dots\dcup V_s$) if its ends lie in different classes of this
partition. A subgraph~$H$ is called crossing if all of its edges are
crossing, and \emph{non-crossing} if none of its edges is crossing.
The \emph{chromatic excess} $\sigma(H)$ of~$H$ denotes the smallest
size of a colour class in a proper $\chi(H)$-colouring of~$H$.

\subsection{Auxiliary lemmas}

The proof of Theorem~\ref{thm:main} relies on the following three lemmas.
Firstly, we use a stability-type result which was observed in~\cite{OP11}.

\begin{lemma}[stability lemma~\cite{OP11}]
\label{lem:stab}
  For every~$\gamma>0$, every integer~$r\ge3$, and every graph $H\neq K_r$ with
  $\chi(H)=r$ there is an integer~$n_0$ such that the following holds. If
  $G=(V,E)$ has $n\ge n_0$ vertices and satisfies $\phi_H(G)\ge\ex(n,K_r)$,
  then there is a partition $V=V_1\dcup\cdots\dcup V_{r-1}$ such that
  \begin{enumerate}[label=\abc]
    \item\label{lem:stab:0} $\deg(v,V_i)\le\deg(v,V_j)$ for all $v\in V_i$
      and all $i,j\in[r-1]$,
    \item\label{lem:stab:1} $\sum_i e(V_i)<\gamma n^2$, and
    \item\label{lem:stab:2} $\frac{n}{r-1}-2\sg n \le |V_i|\le\frac{n}{r-1}+2r\sg n$.
      \qed
  \end{enumerate}
\end{lemma}

We remark that this lemma is stated in~\cite{OP11} only with
assertion~\ref{lem:stab:1}. However, we can certainly assume that the
partition obtained is a maximal $(r-1)$-cut, which implies~\ref{lem:stab:0},
and for \ref{lem:stab:2} see Claim~$8$ in~\cite{OP11}.

The following lemma allows us to find many $H$-copies in a graph~$G$ with a
partition such that each vertex has few neighbours inside its own partition class.

\begin{lemma}
\label{lem:lowdegree}
  For every integer $r\ge 3$, every graph~$H$ with $\chi(H)=r$ and every positive 
  $\beta\le 1/\big(100 e(H)^4\big)$ there is an integer~$n_0$ 
  such that the following holds. Let~$G=(V,E)$ be a graph on $n\ge n_0$
  vertices, with a partition $V=V_1\dcup\dots\dcup V_{r-1}$ such that for
  all $i,j\in[r-1]$ with $i\neq j$
  \begin{enumerate}[label=\rom]
    \item\label{lem:lowdegree:1} $\deg(v,V_j)\ge
      \left(\tfrac{1}{r-1}-\beta\right)n$ for every $v\in V_i$,
    \item\label{lem:lowdegree:2} $\sum_{i'=1}^{r-1} e(V_{i'}) \le\beta^2n^2/e(H)$
      and $\Delta(V_i)\le 2 \beta n$.
\end{enumerate}
Then we can consecutively delete edge-disjoint copies of $H$ from $G$,
until $e(V_i)\le\biex(n,H)$ for all $i\in[r-1]$.  Moreover, these
$H$-copies can be chosen such that each of them contains a non-crossing
$F\in\cFs_H$ and all edges in $E(H)\setminus E(F)$ are crossing.
\end{lemma}


\begin{proof} 
  Let $G=(V,E)$ be a graph and $V=V_1\dcup\dots\dcup V_{r-1}$ be a
  partition satisfying the conditions of the lemma.  We proceed by
  selecting copies of $H$ in~$G$ and deleting them, one at a time, in the
  following way. First we find a copy of some $F\in\cFs_H$ in~$G[V_i]$ for some partition
  class~$V_i$. Then we extend this~$F$ to a copy of~$H$, using only
  vertices~$v$ of~$G$ for $H\setminus F$ which have at least
  $\left(\tfrac{1}{r-1}-2\beta\right)n$ neighbours in every partition class
  other than their own. We say that such vertices~$v$ are
  \emph{$\beta$-active}. 

  We need to show that this deletion process can be
  performed until $e(V_i)\le\biex(n,H)$ for all $i\in[r-1]$. Clearly, while
  $e(V_i)>\biex(n,H)$ for some~$i$, we find some $F\in\cFs_H$ in
  $G[V_i]$. Let such a copy of~$F$ be fixed in the following and assume
  without loss of generality that $V(F)\subset V_{r-1}$. It remains to
  show that~$F$ can be extended to a copy of~$H$.

  By condition~\ref{lem:lowdegree:1}, at the beginning of the deletion
  process every vertex is $\beta$-active, and every vertex which gets
  inactive has lost at least $\beta n$ neighbours in some
  partition class other than its own. Further, by
  condition~\ref{lem:lowdegree:2} we can
  find at most $\beta^2n^2/e(H)$ copies of~$H$ in this way.
  Hence we conclude
  that even after the very last deletion step, the number of vertices which
  are not $\beta$-active is at most 
  \[
    \frac{\beta^2 n^2}{e(H)}e(H)\cdot\frac{1}{\beta n}=\beta n \,.
  \]
  In addition, by condition~\ref{lem:lowdegree:2} we have
  $\Delta(V_{r-1})\le 2\beta n$ at the beginning of the
  deletion process. Recall moreover that, in this process, we 
  use inactive vertices only in copies of some graph in $\cFs_H$ (and not to complete
  such a copy to an $H$-copy). Hence, throughout the process, we have for
  all $j\in[r-1]$ and all $v\in V\setminus V_j$ that
  \begin{equation}\label{eq:lowdegree:deg}
    \deg(v,V_j)
    \ge\big(\tfrac{1}{r-1}-2\beta\big)n-e(H)-e(H)\cdot 2\beta n
    \ge\tfrac{n}{r-1}-5e(H)\beta n\,.
  \end{equation}
  By condition~\ref{lem:lowdegree:1} each partition class~$V_j$ has size at
  least $(\frac1{r-1}-2\beta)n$, and thus size at most
  $\frac{n}{r-1}+2r\beta n$.
  Moreover, by~\eqref{eq:lowdegree:deg} each vertex~$v\in V_j$ has at most
  \begin{equation*}
    \tfrac{n}{r-1}+2r\beta n-\big(\tfrac{n}{r-1}-6e(H)\beta n\big)
    \le 8 e(H) \beta n
  \end{equation*}
  non-neighbours in each~$V_{j'}$ with $j\neq j'$.  Hence, any set
  $S\subset V\setminus V_j$ with $|S|\le r\cdot v(H)$ has at least
  $|V_j|-8r\cdot v(H)e(H)\beta n$ common neighbours in $V_j$. In
  particular, $S$ has at least
  \[(\tfrac1{r-1}-2\beta)n-\beta
  n-8r\cdot v(H)e(H)\beta n\ge \tfrac{n}{r-1}-11e(H)^3\beta n>\beta n\ge v(H)\] common
  neighbours in $V_j$ which are $\beta$-active, where we used the condition
  $\beta\le 1/(100 e(H)^4)$ in the second inequality, and in the last
  inequality that $n$ is sufficiently large.
  
  When~$F$ gets selected in the deletion process, we use the above
  observation to construct within the $\beta$-active common neighbours of
  $F$ a copy of the complete $(r-2)$-partite graph with $v(H)$ vertices in
  each part, as follows.  We inductively find sets $S_i\subset V_i$ of size
  $v(H)$ which form the parts of this complete $(r-2)$-partite graph.  For
  each $1\le i\le r-2$ in turn, we note that $v(F)+(i-1)v(H)\le r\cdot
  v(H)$, and therefore the set $v(F)\cup S_1\cup\cdots\cup S_{i-1}$ has at
  least $\beta n\ge v(H)$ common neighbours in $V_i$ which are
  $\beta$-active.  We let $S_i$ be any set of $v(H)$ of these
  $\beta$-active common neighbours.
  Thus we can extend $F$ to a copy of $H$ in $G$.
\end{proof}

With the help of the next lemma we will find $H$-copies using those
vertices which have many neighbours in their own partition class.

\begin{lemma}\label{lem:highdegree}
  For every integer $r\ge 3$, every 
  graph~$H$ with $\chi(H)=r$, and every positive $\beta\le
  1/\big(2e(H)^2\big)$ there are integers~$K$ and~$n_0$ such that the
  following holds. Let~$G=(V,E)$ be a graph on $n\ge n_0$ vertices, with a
  partition $V=X\dcup V'_1\dcup\dots\dcup V'_{r-1}$
  such that
  \begin{enumerate}[label=\rom]
   \item\label{lem:highdegree:1} $e(V'_i,V'_j)>|V'_i||V'_j|-\beta^6 n^2$ for each $i,j\in[r-1]$ with
      $i\neq j$,
   \item\label{lem:highdegree:2} $|X|\le \beta^6n$. 
  \end{enumerate}
  Then we can consecutively delete edge-disjoint copies of $H$ from $G$,
  until for all but at most $K\big(\sigma(H)-1\big)$ vertices $x \in X$
  there is an $i\in[r-1]$ such that $\deg(x,V'_i)\le\beta^2 n$.  Moreover,
  these $H$-copies can be chosen such that they are crossing for
  the partition $X\dcup V'_1\dcup\dots\dcup V'_{r-1}$ and each of them uses exactly
  exactly $\sigma(H)$ vertices of~$X$.
\end{lemma}
\begin{proof} 
  Without loss of generality we assume that there are only crossing edges
  in~$G$ (otherwise delete the non-crossing edges).  We proceed as
  follows. In the beginning we set~$X':=X$. Then we identify $\sigma(H)$
  vertices in~$X'$ which are completely joined to a complete
  $(r-1)$-partite graph $K_{r-1}\big(v(H)\big)$ in $V\setminus X$, with
  $v(H)$ vertices in each part. The subgraph of~$G$ identified in this way
  clearly contains a copy of~$H$ with the desired properties, whose edges
  we delete from~$G$. Next we delete those vertices~$x$ from~$X'$ with
  $\deg(x,V'_i)\le\beta^2 n$ for some $i\in[r-1]$. Then we continue with
  the next copy of~$H$. 

  We need to show that this process can be repeated
  until $X'\le K\big(\sigma(H)-1\big)$.
  Indeed, assume that we still have $X'>K\big(\sigma(H)-1\big)$.  Observe
  that since $\sum_{x\in X}\deg(x)<|X|n$ we can find less than
  $|X|n\le\beta^6n^2$ copies of~$H$ with the desired properties in total,
  where we used condition~\ref{lem:highdegree:2}.  Hence, throughout the
  process at most $e(H)\beta^6n^2$ edges are deleted from~$G$.  In
  addition, for each $x\in X'$ we have by definition
  $\deg(x,V'_i)>\beta^2 n$ for all $i\in[r-1]$. Hence we can choose for
  each $i$ a set $S_i\subseteq N_{V'_i}(x)$ of size $\beta^2 n$.  By
  condition~\ref{lem:highdegree:1} the graph $G[\dcup S_i]$ has density at
  least
  \begin{equation*}
    \frac{\binom{r-1}{2}
      \big((\beta^2n)^2-\beta^6 n^2-e(H)\beta^6n^2\big)}{\binom{(r-1)\beta^2n}{2}}
    \ge\frac{r-2}{r-1}\big(1-(e(H)+1)\beta^2\big)>\frac{r-3}{r-2} \,,
  \end{equation*}
  where we used $\beta\le 1/\big(2e(H)^2\big)$ in the last inequality.
  Thus, since~$n$ is sufficiently large, we can apply the supersaturation
  theorem of Erd\H{o}s and Simonovits~\cite{ErdSim83}, to conclude that the
  graph $G[\dcup S_i]$ contains at least $\delta n^{(r-1)v(H)}$ copies of
  $K_{r-1}\big(v(H)\big)$, where $\delta>0$ depends only on $\beta$ and
  $e(H)$. Choosing $K:=1/\delta$, we can then use the pigeonhole principle
  and the fact that $|X'|>K\left(\sigma(H)-1\right)$) to infer that there
  are $\sigma(H)$ vertices in~$X'$ which are all adjacent to the vertices
  of one specific copy of $K_{r-1}\big(v(H)\big)$ in $G[\dcup S_i]$ as
  desired.
\end{proof}

In addition we shall use the following easy fact about $\biex(n,H)$.

\begin{fact}\label{fc:biex}
  Let $H$ be an $r$-chromatic graph, $r\ge 3$. 
  If $\biex(n,H)<n-1$ then $\sigma(H)=1$.
\end{fact}
\begin{proof}
  If $\sigma(H)\ge 2$, then each $F$ from $\cFs_H$ contains a matching of
  size~$2$. Thus $\biex(n,H)\ge n-1$ since the star $K_{1,n-1}$ does
  not contain two disjoint edges. 
\end{proof}

\section{Proof of Theorem~\ref{thm:main}}
\label{sec:proof}

In this section we show how Lemmas~\ref{lem:stab}, \ref{lem:lowdegree} and
\ref{lem:highdegree} imply Theorem~\ref{thm:main}.

\begin{proof}[Proof of Theorem~\ref{thm:main}]
  Let~$r$ and~$H$ with $\chi(H)=r\ge 3$ be given.
  If~$H=K_r$, then the result of Bollob\'as~\cite{Bo76}
  applies, hence we can assume that $H\neq K_r$. We choose
  \begin{equation}\label{eq:main:betagamma}
    \beta:= \frac{1}{100e(H)^4} \qqand
    \gamma:= \frac{\beta^{12}}{1000 e(H)^4} \,.
  \end{equation}
  Let~$K$ be the constant from Lemma~\ref{lem:highdegree} and choose
  \begin{equation}\label{eq:main:C}
    C:=K \cdot v(H) \beta^{-1} \,.
  \end{equation}
  Finally let~$n_0$ be sufficiently large for Lemmas~\ref{lem:stab},
  \ref{lem:lowdegree} and~\ref{lem:highdegree}.

  Now let $G$ be a graph with $n\ge n_0$ vertices and 
  \begin{equation}\label{eq:main:eG}
    e(G)\ge\ex(n,K_r)+C\cdot\biex(n,H)\,,
  \end{equation}
  and assume for contradiction that 
  \begin{equation}\label{eq:main:phi}
    \phi_H(G)\ge\ex(n,K_r)\,.
  \end{equation}

  Observe first that we may assume without loss of generality that 
  \begin{equation}\label{eq:main:delta}
    \delta(G)\ge\delta\big(T_{r-1}(n)\big)\,.
  \end{equation}
  Indeed, if this is not the case, we can consecutively delete vertices of
  minimum degree until we arrive at a graph $G_{n^*}$ on~$n^*$ vertices
  with $\delta(G_{n^*})\ge \delta\big(T_{r-1}(n^*)\big)$. Denote the
  sequence of graphs obtained in this way by $G_n:=G$, $G_{n-1}$, \ldots,
  $G_{n^*}$.  We have
  \[
    \ex(n,K_r)\le\phi_H(G)\le \phi_H(G_{n-1})+\delta\big(T_{r-1}(n)\big)-1
  \]
  and thus $\phi_H(G_{n-1})\ge\ex(n-1,K_r)+1$. Similarly $\phi_H(G_{n-i})\ge
  \ex(n-i,K_r)+i$.  Since~$n$ is sufficiently large there is an~$i^*$ such
  that $n-i^*\ge n_0$ and $i^*\ge\binom{n-i^*}{2}+1$. Hence $n^*>n-i^*\ge
  n_0$, since otherwise $\phi_H(G_{n^*})\ge \ex(n^*,K_r)+\binom{n^*}{2}+1$,
  a contradiction. Thus we may assume~\eqref{eq:main:delta}.

  Next, by~\eqref{eq:main:phi}, we can apply Lemma~\ref{lem:stab}, which
  provides us with a partition $V_1\dcup\ldots\dcup
  V_{r-1}$ of $V(G)$ such that assertions~\ref{lem:stab:0},~\ref{lem:stab:1} and~\ref{lem:stab:2} 
  in Lemma~\ref{lem:stab} are satisfied.  
  Let $m:=\sum_{i=1}^{r-1} e(V_i)$.
  Equation~\eqref{eq:main:eG} and Lemma~\ref{lem:stab}\ref{lem:stab:1} imply
  \begin{equation}\label{eq:main:m}
     C\cdot\biex(n,H) \le m\le \gamma n^2 
     \leByRef{eq:main:betagamma} \beta^2 n^2 / e(H) \,.
  \end{equation}
  Further, by the definition of~$m$ we clearly have
  $e(G)\le\ex(n,K_r)+m$.
  Hence it will suffice to find
  $\frac{m}{e(H)-1}+1$
  edge-disjoint copies of~$H$ in~$G$, since this would imply
  \begin{equation*}\label{eq:contradiction}
    \phi_H(G)\le \ex(n,K_r)+m -\left(\frac{m}{e(H)-1}+1\right)\big(e(H)-1\big)<\ex(n,K_r)\,,
  \end{equation*}
  contradicting~\eqref{eq:main:phi}. So this will be our goal in the
  following, which we shall achieve by first applying
  Lemma~\ref{lem:lowdegree} and then Lemma~\ref{lem:highdegree}.

  We prepare these applications by identifying
  for every $i\in[r-1]$ the set~$X_i$ of vertices in~$V_i$ with high
  degree to its own class, that is,
  \[
  X_i:=\big\{v\in V_i\colon \deg(v,V_i)\ge \tfrac12\beta n\big\}\,.
  \]
  Let~$X:=\dcup_{i\in[r-1]}X_i$.
  This implies
 \begin{equation}\label{eq:main:X}
    |X|\le \frac{2m}{\frac12 \beta n}
    \leByRef{eq:main:m} \frac{2\gamma n^ 2}{\frac12 \beta n}
    \leByRef{eq:main:betagamma} \sqrt\gamma n
    \leByRef{eq:main:betagamma} \beta^6 n\,.
 \end{equation}
 In addition we set $V'_i:= V_i\setminus X_i$ for all $i\in[r-1]$,
 $n':=|V\setminus X|$, $m':=\sum_{i=1}^{r-1} e(V_i\setminus X_i)$ and
 $m_X:=m-m'=e(X)+\sum_{i=1}^{r-1}e(X_i,V'_i)$. 

 \smallskip

\noindent \textbf{Step~1.}  We want to apply
 Lemma~\ref{lem:lowdegree} to the graph $G[V\setminus X]$ and the partition
 $V'_1\dcup\dots\dcup V'_{r-1}$. We first need to
 check that the conditions are satisfied. By Lemma~\ref{lem:stab}\ref{lem:stab:2}
 and~\eqref{eq:main:X} we have for each for each $i,j\in[r-1]$
 with $i\neq j$ that $|V\setminus(V_i\cup
 V'_j)|\le(\frac{r-3}{r-1}+6\sqrt{\gamma})n$.
 Moreover, by~\eqref{eq:main:X} we clearly have $n'\ge n/2$.
 Hence, by the definition of~$X$, for each $v\in V'_i$ we have
  \begin{equation}\label{eq:main:deg}
    \begin{split}
    \deg(v,V'_j)
    & \geByRef{eq:main:delta} \delta\big(T_{r-1}(n)\big) -
    |V\setminus(V_i\cup V'_j|)| - \tfrac12\beta n \\
    & \ge \big(\tfrac{1}{r-1}-7\sqrt{\gamma}-\tfrac12\beta\big) n
    \geByRef{eq:main:betagamma} \big(\tfrac{1}{r-1}-\beta\big) n'
    \end{split}
 \end{equation}
  and thus condition~\ref{lem:lowdegree:1} of Lemma~\ref{lem:lowdegree} is satisfied. 
  Condition~\ref{lem:lowdegree:2} of Lemma~\ref{lem:lowdegree} holds
  by~\eqref{eq:main:m} and the definition of~$X$.
  Therefore we can apply Lemma~\ref{lem:lowdegree}.

  This lemma asserts that we can consecutively delete copies of~$H$ from
  $G[V\setminus X]$, each containing a non-crossing $F\in\cFs_H$ and
  crossing edges otherwise, until $e(V'_i)\le\biex(n',H)$. Denote the graph
  obtained after these deletions by~$G_1$.
 
  We have $\max_{F\in\cFs_H} e(F)\le e(H)-2$, since~$\chi(H)\ge3$.  Hence each
  copy of~$H$ deleted in this way uses at most $e(H)-2$ non-crossing edges, and so this gives at least
  \begin{equation}\label{eq:firststep}
    \frac{m'-(r-1)\biex(n,H)}{e(H)-2}\ge\frac{m'-r\biex(n,H)}{e(H)-2}
  \end{equation}
  edge-disjoint copies of~$H$ in $G[V\setminus X]$.  
  
  By assertion~\ref{lem:stab:1} of Lemma~\ref{lem:stab} and the assumption
  $e(G)\ge\ex(n,K_r)$, we have that $e_G(V_i,V_j)\ge |V_i||V_j|-\gamma
  n^2$, and thus $e_G(V'_i,V'_j)\ge|V'_i||V'_j|-\gamma n^2$.  Again by
  assertion~\ref{lem:stab:1} of Lemma~\ref{lem:stab}, in obtaining $G_1$ as
  described above we delete at most $\gamma n^2$ copies of $H$, so we have
  \begin{equation}\label{eq:densepairs}
 \begin{split}
   e_{G_1}(V'_i,V'_j)
   & \ge 
   |V'_i||V'_j|-\gamma n^2-(e(H)-1)\gamma n^2 \\
   &=|V'_i||V'_j|-e(H)\gamma n^2   
    \geByRef{eq:main:betagamma} |V'_i||V'_j| - \beta^6 n^2 \,.
  \end{split}
  \end{equation}

 \smallskip\noindent \textbf{Step~2.} Next we want to apply
  Lemma~\ref{lem:highdegree} to ~$G_1$ and the partition $X\dcup V'_1\dcup\dots\dcup
  V'_{r-1}$. Note that condition~\ref{lem:highdegree:1} of
  Lemma~\ref{lem:highdegree} is satisfied by~\eqref{eq:densepairs} and
  condition~\ref{lem:highdegree:2} by~\eqref{eq:main:X}. Hence
  Lemma~\ref{lem:highdegree} allows us to delete crossing copies of $H$
  from~$G$ until all vertices~$x$ of a subset $X_0\subset X$ with
  $|X|-|X_0|\le K(\sigma(H)-1)=:K'$ have $\deg(x,V'_{i(x)})\le\beta^2 n$ for some
  $i(x)\in[r-1]$. Denote the graph obtained after these deletions by~$G_2$.

  Now let $x\in X_j$ for some $j\in[r-1]$ be arbitrary. We set
  $m_x:=\deg_{G}(x,V_j\setminus X)$.  Since no edges adjacent to~$x$ were
  deleted in step~1, if $x\in X_0$ then the number of edges adjacent to~$x$
  deleted in step~2 is at least $\deg_{G}(x,V_{i(x)}\setminus X)-\beta^2
  n\ge m_x-2\beta^2 n$, where we used assertion~\ref{lem:stab:0} of
  Lemma~\ref{lem:stab} and~\eqref{eq:main:X} in the inequality. Hence,
  since $m_X=\sum_{x\in X}m_x + e(X)$, in total at least
 \begin{multline*} 
  m_X -K'n
  -|X_0|\beta^2 n
  -e(X)
  \ge m_X-K'n-2\beta^2n|X| \\
  \geByRef{eq:main:X} m_X-K'n-8\beta m
 \end{multline*}
  edges adjacent to~$X$ were deleted in step~$2$.
  By Fact~\ref{fc:biex} we have $K'=K(\sigma(H)-1)=0$ if
  $\biex(n,H)<n-1$. If $\biex\ge n-1\ge n/2$ on the other hand, then
  $m\ge Cn/2$ by~\eqref{eq:main:m} and thus
  $K' n\le 2K' m/ C $.
  Observe moreover that, because~$H\neq K_3$, each
  $H$-copy deleted in this step uses at least~$2$ edges which are not adjacent to~$X$.
  We conclude that at least
 \begin{equation}\label{eq:secstep}
      \frac{m_X-\frac{2K(\sigma(H)-1)}{C}m -8\beta m}{e(H)-2}
      \geByRef{eq:main:C}\frac{m_X-9\beta m}{e(H)-2}
  \end{equation}
  edge-disjoint copies of $H$ were deleted from~$G_1$ in step~$2$. 

  \smallskip

  Combining~\eqref{eq:firststep} and~\eqref{eq:secstep}
  reveals that~$G$ contains
  \begin{multline*}
   \frac{m'-r\biex(n,H)}{e(H)-2} + \frac{m_X-9\beta m}{e(H)-2} 
   \geByRef{eq:main:m} \frac{m-\frac{r}{C}m-9\beta m}{e(H)-2} \\ 
   \geByRef{eq:main:C} \frac{m-10\beta m}{e(H)-2}
   \geByRef{eq:main:betagamma} \frac{m}{e(H)-1}+1
  \end{multline*}
  edge-disjoint copies of $H$, which gives the desired contradiction.
\end{proof}


\bibliographystyle{amsplain}
\bibliography{biexphi}

\end{document}